\documentclass[11pt]{amsart}
\usepackage{amsmath, amsfonts, eucal, amsbsy, amssymb, upref, graphicx, chatterjee1}

\newtheorem{thm1}{Theorem}

\begin{document}
\title{An observation about submatrices}
\author{Sourav Chatterjee  }
\address{\newline367 Evans Hall \#3860\newline
Department of Statistics\newline
University of California at Berkeley\newline
Berkeley, CA 94720-3860\newline
{\it E-mail: \tt sourav@stat.berkeley.edu}\newline 
}
\thanks{Sourav Chatterjee's research was partially supported by NSF grant DMS-0707054 and a Sloan Research Fellowship}
\author{Michel Ledoux}
\address{\hskip-.15in Institut de Math\'ematiques \newline
Universit\'e de Toulouse \newline
31062 Toulouse Cedex 9 France\newline
{\it E-mail: \tt ledoux@math.univ-toulouse.fr}
}
\thanks{Michel Ledoux's research was partially supported by the ANR Grandes Matrices Al\'eatoires}
\keywords{Random matrix, concentration of measure, empirical distribution, eigenvalue}
\subjclass[2000]{60E15, 15A52}
\begin{abstract}
Let $M$ be an arbitrary Hermitian matrix of order $n$, and $k$ be a positive integer $\le n$. We show that if $k$ is large, the distribution of eigenvalues on the real line is almost the same for almost all  principal submatrices of $M$ of order $k$. 
The proof uses results about random walks  on symmetric groups and concentration of measure.
In a similar way, we also show that almost all $k\times n$ submatrices of $M$ have almost the same distribution of singular values.
\end{abstract}
\maketitle 

Let $M$ be a square matrix of order $n$. For any two sets of integers $i_1,\ldots, i_k$ and $j_1,\ldots, j_l$ between $1$ and $n$, $M(i_1,\ldots,i_k; j_1,\ldots,j_l)$ denotes the submatrix of $M$ formed by deleting all rows except rows $i_1,\ldots, i_k$, and all columns except columns $j_1,\ldots, j_l$. A submatrix like $M(i_1,\ldots,i_k; i_1,\ldots,i_k)$ is called a principal submatrix. 

For a Hermitian matrix $M$ of order $n$ with eigenvalues $\lambda_1,\ldots, \lambda_n$ (repeated by multiplicities), let $F_M$ denote the empirical spectral distribution function of $M$, that is, 
\[
F_M(x) := \frac{\#\{i: \lambda_i \le x\}}{n}.
\]
The following result shows that given $1\ll k\le n$ and any Hermitian matrix $M$ of order $n$, the empirical spectral distribution is {\it almost the same} for {\it almost every} principal submatrix of $M$ of order $k$. 
\begin{thm1}\label{mainthm}
Take any $1\le k\le n$ and a Hermitian matrix $M$ of order $n$. Let $A$ be a principal submatrix of $M$ chosen uniformly at random from the set of all $k\times k$ principal submatrices of $M$. Let $F$ be the expected spectral distribution function of $A$, that is, $F(x) = \ee F_A(x)$. Then for each $r\ge 0$,
\[
\pp(\|F_A - F\|_\infty \ge k^{-1/2} + r) \le 12\sqrt{k}e^{-r\sqrt{k/8}}.
\] 
Consequently, we have
\begin{align*}
\ee\|F_A - F\|_\infty &\le  \frac{13 + \sqrt{8}\log k}{\sqrt{k}}.
\end{align*}
Exactly the same results hold if $A$ is a $k\times n$ submatrix of $M$ chosen uniformly at random, and $F_A$ is the empirical distribution function of the singular values of $A$. Moreover, in this case $M$ need not be Hermitian.
\end{thm1}
\noindent {\it Remarks.} (i) Note that the bounds do not depend at all on the entries of~$M$, nor on the dimension $n$.

(ii) We think it is possible to improve the $\log k$ to $\sqrt{\log k}$ using Theorem~2.1 of Bobkov~\cite{bobkov04} instead of the spectral gap techniques that we use. (See also Bobkov and Tetali \cite{bobkovtetali06}.) However, we do not attempt to make this small improvement because $\sqrt{\log k}$, too, is unlikely to be optimal. Taking $M$ to be the matrix which has $n/2$ $1$'s on the diagonal and the rest of the elements are zero, it is easy to see that there is a lower bound of $const. k^{-1/2}$. We conjecture that the matching upper bound is also true, that is, there is a universal constant $C$ such that $\ee\|F_A - F\|_\infty \le Ck^{-1/2}$.

(iii) The function $F$ is determined by $M$ and $k$.  If $M$ is a diagonal matrix, then $F$ is exactly equal to the spectral measure of $M$, irrespective of~$k$. However it is not difficult to see that the spectral measure of $M$ cannot, in general, be reconstructed from $F$. 

(iv) The result about random $k\times n$ submatrices is related to the recent work of Rudelson and Vershynin \cite{rudelsonvershynin07}. Let us also refer to \cite{rudelsonvershynin07} for an extensive list of references to the substantial volume of  literature on random submatrices in the computing community. However, most of this literature (and also \cite{rudelsonvershynin07}) is concerned with the largest eigenvalue and not the bulk spectrum. On the other hand, the existing techniques are usually applicable  only when $M$ has low rank or low `effective rank' (meaning that most eigenvalues are negligible compared to the largest one).
\vskip.1in

\begin{figure}
\centering \begin{tabular}{c} \resizebox{!}{7cm}{\includegraphics{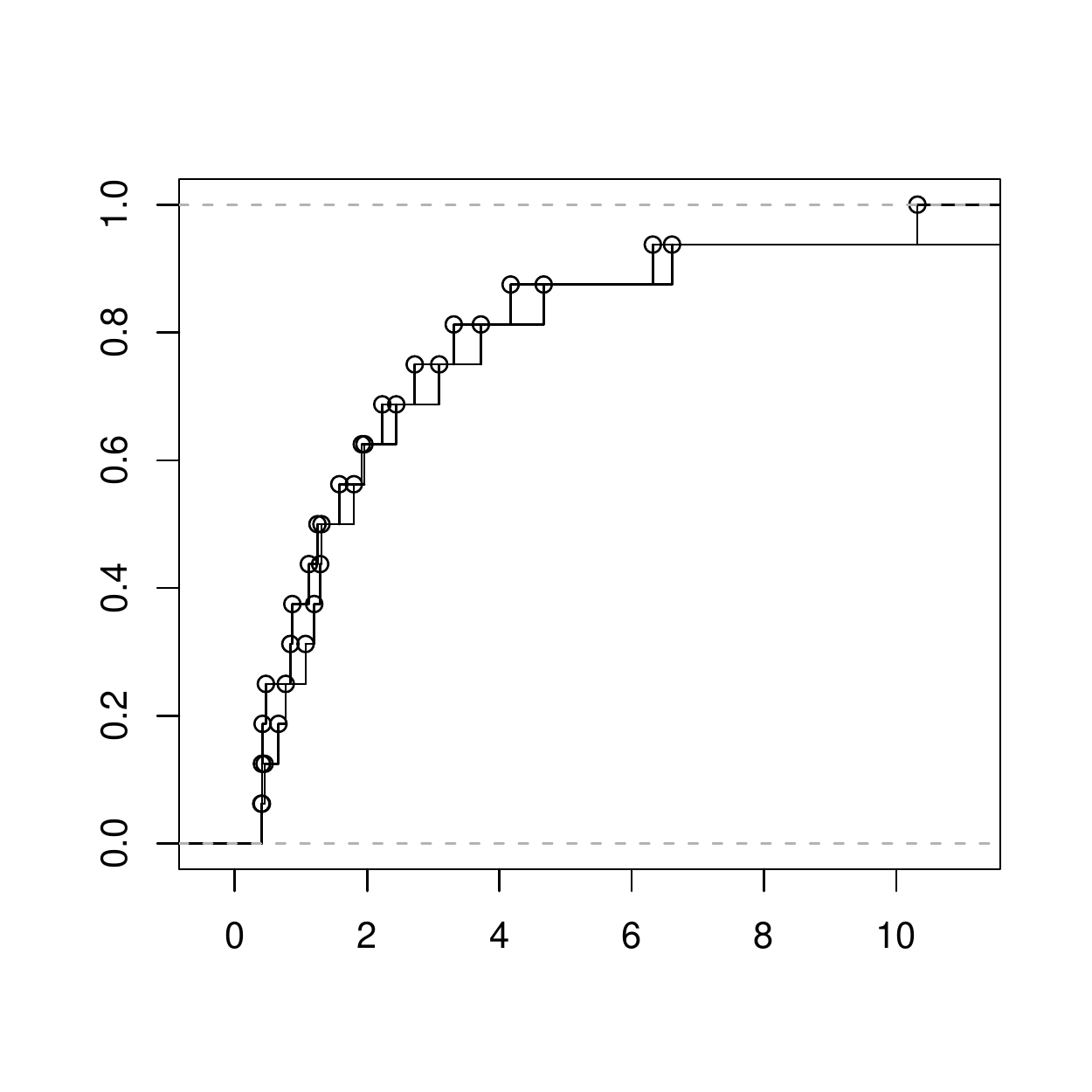}} 
\end{tabular}
\caption{Superimposed empirical distribution functions of two submatrices of order $20$ chosen at random from a deterministic matrix of order $100$.} 
  \label{fig1}
\end{figure}

{\it A numerical illustration.} The following simple example demonstrates that the effects of Theorem \ref{mainthm} can kick in even when $k$ is quite small. We took $M$ to be a $n\times n$ matrix for $n = 100$, with $(i,j)$th entry = $\min\{i,j\}$. This is  the covariance matrix of a simple random walk up to time $n$. We chose $k = 20$, and picked two $k\times k$ principal submatrices $A$ and $B$ of $M$, uniformly and independently at random. Figure \ref{fig1} plots to superimposed empirical distribution functions of $A$ and $B$, after excluding the top $4$ eigenvalues since they are too large. The classical Kolmogorov-Smirnov test from statistics gives a $p$-value of $0.9999$ (and $\|F_A - F_B\|_\infty = 0.1$), indicating that the two distributions are {\it statistically indistinguishable}.

\vskip.1in

{\it Markov chains.} Let us now quote two results about Markov chains that we need to prove Theorem \ref{mainthm}.
Let $\mx$ be a finite or countable set. Let $\Pi(x,y)\ge 0$ satisfy
\[
\sum_{y\in \mx} \Pi(x,y) = 1
\]
for every $x\in \mx$. Assume furthermore that there is a symmetric invariant probability measure $\mu$ on $\mx$, that is, $\Pi(x,y)\mu(\{x\})$ is symmetric in $x$ and $y$, and $\sum_x \Pi(x,y) \mu(\{x\}) = \mu(\{y\})$ for every $y\in \mx$. In other words, $(\Pi,\mu)$ is a reversible Markov chain. For every $f:\mx \ra \rr$, define
\[
\mathcal{E}(f,f) = \frac{1}{2}\sum_{x,y\in \mx} (f(x)-f(y))^2 \Pi(x,y) \mu(\{x\}).
\]
The spectral gap or the Poincar\'e constant of the chain $(\Pi,\mu)$ is the largest $\lambda_1 > 0$ such that for all $f$'s, 
\[
\lambda_1 \var_\mu(f) \le \mathcal{E}(f,f).
\]
Set also
\begin{equation}\label{trinorm}
||| f |||_\infty^2 = \frac{1}{2}\sup_{x\in \mx}\sum_{y \in \mx} (f(x)-f(y))^2 \Pi(x,y).
\end{equation}
The following concentration result is a copy of Theorem 3.3 in \cite{ledoux01}.
\begin{thm1}[\cite{ledoux01}, Theorem 3.3]\label{ledoux}
Let $(\Pi,\mu)$ be a reversible Markov chain on a finite or countable space $\mx$ with a spectral gap $\lambda_1 > 0$. Then, whenever $f:\mx \ra \rr$ is a function such that  $|||f|||_\infty \le 1$, we have that $f$ is integrable with respect to $\mu$ and for every $r\ge 0$,
\[
\mu(\{f\ge {\textstyle \int} f d\mu + r\}) \le 3 e^{-r\sqrt{\lambda_1}/2}.
\]
\end{thm1}
\noindent Let us now specialize to $\mx = S_n$, the group of all permutations of $n$ elements. The following transition kernel $\Pi$ generates the `random transpositions walk'.
\begin{equation}\label{pidef}
\Pi(\pi,\pi') =
\begin{cases}
1/n &\text{ if } \pi' = \pi,\\
2/n^2 &\text{ if } \pi' = \pi \tau \text{ for some transposition } \tau,\\
0 &\text{ otherwise.}
\end{cases}
\end{equation}
It is not difficult to verify that the uniform distribution $\mu$ on $S_n$ is the unique invariant measure for this kernel, and the pair $(\Pi,\mu)$ defines a reversible Markov chain. 
\begin{thm1}
[Diaconis \& Shahshahani \cite{diaconisshahshahani81}, Corollary 4]\label{ds}
The spectral gap of the random transpositions walk on $S_n$ is $2/n$. 
\end{thm1}
\vskip.1in
We are now ready to prove Theorem \ref{mainthm}.

\begin{proof}[Proof of Theorem \ref{mainthm}]
Let $\pi$ be a uniform random permutation of $\{1,\ldots,n\}$. Let $A = A(\pi) = M(\pi_1,\ldots,\pi_k; \pi_1,\ldots,\pi_k)$. Fix a point $x\in \rr$. Let 
\[
f(\pi) := F_{A}(x).
\]
Let $\Pi$ be the transition kernel for the random transpositions walk defined in \eqref{pidef}, and let $||| \cdot|||_\infty$ be defined as in \eqref{trinorm}.

Now, by Lemma 2.2 in Bai \cite{bai99}, we know that for any two Hermitian matrices $A$ and $B$ of order $k$,
\begin{equation}\label{s1}
\|F_A - F_B\|_\infty \le \frac{\mathrm{rank}(A-B)}{k}.
\end{equation}
Let $\tau = (I,J)$ be a random transposition, where $I,J$ are chosen independently and uniformly from $\{1,\ldots,n\}$. Multiplication by $\tau$ results in taking a step in the chain defined by $\Pi$. Now, for any $\sigma\in S_n$, the $k\times k$  Hermitian matrices $A(\sigma)$ and $A(\sigma \tau)$ differ at most in one column and one row, and hence $\mathrm{rank}(A(\sigma)-A(\sigma\tau))\le 2$. Thus,
\begin{equation}\label{s2}
|f(\sigma) - f(\sigma\tau)|\le \frac{2}{k}.
\end{equation}
Again, if $I$ and $J$ both fall outside $\{1,\ldots,k\}$, then $A(\sigma) = A(\sigma\tau)$. Combining this with \eqref{s1} and \eqref{s2}, we get
\[
|||f|||_\infty^2 = \frac{1}{2}\max_{\sigma\in S_n} \ee(f(\sigma)-f(\sigma\tau))^2 \le \frac{1}{2}\biggl(\frac{2}{k}\biggr)^2 \frac{2k}{n}\le \frac{4}{kn}.
\]
Therefore, from Theorems \ref{ledoux} and \ref{ds}, it follows that for any $r\ge 0$,
\begin{equation}\label{tail1}
\pp(|F_A(x) - F(x)| \ge r) \le 6 \exp\biggl(-\frac{r\sqrt{2/n}}{2\sqrt{4/kn}}\biggr) = 6\exp\biggl(-\frac{r\sqrt{k}}{\sqrt{8}}\biggr).
\end{equation}
The above result is true for any $x$.  Now, if $F_A(x-) := \lim_{y\uparrow x} F_A(y)$, then by the bounded convergence theorem we have $\ee F_A(x-) = \lim_{y\uparrow x} F(y) = F(x-)$. It follows that for every $r$, 
\begin{align*}
\pp(|F_A(x-) - \ee F_A(x-)| > r) &\le \liminf_{y\uparrow x} \pp(|F_A(y) - F(y)| > r) \\
&\le 6\exp\biggl(-\frac{r\sqrt{k}}{\sqrt{8}}\biggr).
\end{align*}
Since this holds for all $r$, the $>$ can be replaced by $\ge$. 
Similarly it is easy to show that $F$ is a legitimate cumulative distribution function. Now fix an integer $l\ge 2$, and for $1\le i < l$ let
\[
t_i := \inf \{x: F(x)\ge i/l\}.
\]
Let $t_0=-\infty$ and $t_l = \infty$. Note that for each $i$, $F(t_{i+1}-) - F(t_i)\le 1/l$. Let
\[
\Delta = (\max_{1\le i < l} |F_A(t_i)-F(t_i)| )\vee (\max_{1\le i< l}|F_A(t_i-)- F(t_i-)|).
\]
Now take any $x\in \rr$. Let $i$ be an index such that  $t_i\le x < t_{i+1}$. Then
\[
F_A(x) \le F_A(t_{i+1}-)\le F(t_{i+1}-) + \Delta \le F(x) + 1/l + \Delta.
\]
Similarly,
\[
F_A(x) \ge F_A(t_i) \ge F(t_i) - \Delta \ge F(x)-1/l - \Delta.
\]
Combining, we see that
\[
\|F_A - F\|_\infty \le 1/l + \Delta.
\]
Thus, for any $r\ge 0$,
\[
\pp(\|F_A - F\|_\infty \ge 1/l + r) \le 12(l-1)e^{-r\sqrt{k/8}}.
\]
Taking $l = [k^{1/2}] + 1$, we get for any $r\ge 0$, 
\[
\pp(\|F_A - F\|_\infty \ge 1/\sqrt{k} + r) \le 12\sqrt{k}e^{-r\sqrt{k/8}}.
\] 
This proves the first claim of Theorem \ref{mainthm}. 
To prove the second, using the above inequality, we get
\begin{align*}
\ee\|F_A - F\|_\infty &\le \frac{1 + \sqrt{8}\log k}{\sqrt{k}} +  \pp\biggl(\|F_A - F\|_\infty \ge \frac{1+ \sqrt{8}\log k}{\sqrt{k}}\biggr)\\
&\le \frac{13 + \sqrt{8}\log k}{\sqrt{k}}.
\end{align*}
For the case of singular values, we proceed as follows. As before, we let $\pi$ be a random permutation of $\{1,\ldots,n\}$; but here we define $A(\pi) = M(\pi_1,\ldots, \pi_k; 1,\ldots, n)$. Since singular values of $A$ are just square roots of eigenvalues of $AA^*$, therefore 
\[
\|F_A - \ee(F_A)\|_\infty = \|F_{AA^*} - \ee(F_{AA^*})\|_\infty, 
\]
and so it suffices to prove a concentration inequality for $F_{AA^*}$. As before, we fix $x$ and define
\[
f(\pi) = F_{AA^*}(x).
\]
The crucial observation is that by Lemma 2.6 of Bai \cite{bai99}, we have that for any two $k \times n$ matrices $A$ and $B$,
\[
\|F_{AA^*} - F_{BB^*}\|_\infty \le \frac{\mathrm{rank}(A-B)}{k}.
\]
The rest of the proof proceeds exactly as before.
\end{proof}

\vskip.1in
\noindent {\bf Acknowledgment.} We thank the referees for helpful comments.

\end{document}